\documentclass[12pt]{article}

\usepackage{amssymb,amsmath,amsfonts,enumerate}

\title{
\author{\bf{Lech Pasicki}}
\bf{Istances and Br\o ndsted's variational principle}
\footnote{MSC: 49J27, 49J45, 54H25, 54E15, 06A06}
\footnote{Author is indebted to professor Christiane Tammer for her valuable suggestions.}}

\newtheorem{theorem}{\indent Theorem}
\newtheorem{lemma}[theorem]{\indent Lemma}
\newtheorem{proposition}[theorem]{\indent Proposition}
\newtheorem{definition}[theorem]{\indent Definition}

\newcommand{\mR}{\mbox{$\mathbb{R}$}}
\newcommand{\mN}{\mbox{$\mathbb{N}$}}
\newcommand{\xnn}{\mbox{$(x_{n})_{n \in \mathbb{N}}$}}
\newcommand{\Xq}{\mbox{$(X,q)$}}
\newcommand{\no}{\mbox{$n_{0}$}}
\newcommand{\qxnxm}{\mbox{$q(x_{n},x_{m})$}}
\newcommand{\XU}{\mbox{$(X,\mathcal{U})$}}
\newcommand{\U}{\mbox{$\mathcal{U}$}}

\newcommand{\C}{\mbox{$\mathcal{C}$}}
\newcommand{\Cy}{\mbox{$C_{y}$}}
\newcommand{\clCy}{\mbox{$\overline{C_{y}}$}}
\newcommand{\xyzX}{\mbox{$x,y,z \in X$}}
\newcommand{\xyX}{\mbox{$x,y \in X$}}
\newcommand{\psiX}{\mbox{$\psi\colon X \to \mathbb{R}$}}
\newcommand{\preq}{\mbox{$\preccurlyeq$}}
\newcommand{\yCx}{\mbox{$y \in C \setminus \{x\}$}}
\newcommand{\yXx}{\mbox{$y \in X \setminus \{x\}$}}
\newcommand{\xo}{\mbox{$x_{0}$}}
\newcommand{\psixnn}{\mbox{$(\psi(x_{n}))_{n \in \mathbb{N}}$}}
\newcommand{\xm}{\mbox{$x_{m}$}}
\newcommand{\qxm}{\mbox{$q(\cdot,x_{m})$}}
\newcommand{\psial}{\mbox{$\psi_{\alpha}$}}
\newcommand{\xn}{\mbox{$x_{n}$}}
\newcommand{\xnp}{\mbox{$x_{n+1}$}}

\newcommand{\FXY}{\mbox{$F\colon X \to 2^{Y}$}}
\newcommand{\xFx}{\mbox{$x \in X \setminus F(x)$}}
\newcommand{\li}{\mbox{$\lim_{n \rightarrow \infty}$}}

\begin{document}
\maketitle
\vspace{1 in}

\begin{abstract}
\par The manuscript contains elegant extensions of the fundamental 
variational principles: Br\o ndsted's and Ekeland's. On the other hand, we get general 
and precise version of the Takahashi and the Caristi fixed point theorems. 
The results are based on the notion of istance for uniform spaces.
\end{abstract}

\par The main results of this paper are theorems of variational type for uniform spaces. Their proofs 
involve a kind of order relation $y \preq x$ defined by $\psi(y) + q(y,x) - \psi(x) \leq 0$. 
Here $q$, called istance, is a special mapping satisfying the triangle inequality (see Definition \ref{De2}).
In consequence we get natural and more general versions of Br\o ndsted's variational principle (Theorem 
\ref{Th6}) and of Ekeland's one (Theorem \ref{Th10}). On the other hand, we obtain a good-looking extension 
of the Caristi and Takahashi fixed point theorems (Theorem \ref{Th8}).
\par Many advanced results for uniform spaces are known. The idea is to replace metric by a family of semimetrics 
(see e.g. \cite{TuM}), pseudometrics (see e.g. \cite{Mi}) or quasimetris (see e.g. \cite{Ha}); and further, 
real-valued mappings by vector-valued functions (see e.g. \cite{HL}, \cite{TuV}). Our theorems look more classical 
while the assumptions on order relation (only transitivity) and $q$ are weak. Theorem \ref{Th6} and in some 
sense similar result of V\'alyi \cite[Theorem 5]{Va} are compared further.
\par Let us recall \cite[Definition 5]{Pv}.
\begin{definition}
\label{De1}
  Let $q\colon X \times X \to [0,\infty)$ be a mapping. Then $\xnn$ is a {\bf Cauchy sequence} 
in $\Xq$, if for  each $\epsilon > 0$ there exists an $\no \in \mN$ such that each 
$m,n \in \mN$, $\no < m < n$ yield $\qxnxm < \epsilon$.
\end{definition}
\par The notion of istance for metric spaces presented in \cite[Definition 6]{Pv} can be extended to the 
case of uniform spaces (see \cite[p. 176]{Ke}) as follows:
\begin{definition}
\label{De2}
  Let $\XU$ be a uniform space. A mapping $q\colon X \times X \to [0,\infty)$ is a 
$\U${\bf-istance} in $X$ if the following conditions are satisfied:
\begin{subequations}
\label{con1}
\begin{align}
 & q(z,x) \leq q(z,y) + q(y,x), \quad \xyzX, \label{con1a}\\
 & q(\cdot,x) \mbox{ is lower semicontinuous,} \quad x \in X, \label{con1b}\\
 & \mbox{each Cauchy sequence in } \Xq \mbox{ is a Cauchy sequence in } \XU. \label{con1c}
\end{align}
\end{subequations} 
\end{definition}
\par The next proposition extends \cite[Proposition 19]{Pv}.
\begin{proposition}
\label{Pro3}
Let $\XU$ be a uniform Hausdorff space. Then any $\U$-istance $q$ in $X$ satisfies
\begin{equation}
\label{con2}
 q(x,y)=0 \mbox{ and } q(y,x)=0 \mbox{ yield } x=y,  \quad \xyX.
\end{equation}
\end{proposition}
\par Proof.
Assume $q(x,y)=0$ and $q(y,x)=0$. Let us adopt $x_{2k-1} = x$ and $x_{2k} = y$, 
$k \in \mN$. Clearly, $\qxnxm = 0$ holds for each $m \neq n$, $m, n \in \mN$, and 
$\xnn$ is a Cauchy sequence in $\XU$ (see (\ref{con1c})), i.e. $x = y$. $\quad \square$
\par If $q\colon X \times X \to \mR$, $\psiX$ are mappings and $q$ satisfies (\ref{con1a}), then 
the following condition defines a transitive relation $\preq$:
\begin{equation}
\label{con3}
 y \preq x \mbox{ iff } \psi(y) + q(y,x) - \psi(x) \leq 0, \quad \xyX.
\end{equation}
\begin{lemma}
\label{Le4}
Let $\XU$ be a complete uniform Hausdorff space, $q$ a $\U$-istance in $X$, and let $\psiX$ 
be a lower semicontinuous mapping bounded below. Assume that $C \subset X$ is a nonempty 
maximal chain (for $\preq$) such that $q(x,y)=0$ or $q(y,x)=0$, for each $x,y \in C$,
$x \neq y$. Then $C$ has a unique smallest element.
\end{lemma}
\par Proof.
Assume that $C$ contains at least two elements. Let us adopt $\Cy=\{z \in C\colon q(z,y)=0 \}$ and suppose 
that there exists a $U \in \U$ such that $\Cy \not \subset U(y)$, for each $y \in C$. Then 
there exists a sequence $\xnn$ in $C$ such that $q(\xnp,\xn)=0$ and $\xnp \not \in U(\xn)$. 
In view of (\ref{con1a}), $\xnn$ is a Cauchy sequence in $\Xq$ and consequently, $\xnn$ is a Cauchy sequence 
in $\XU$. This contradicts the way $\xnn$ was defined. Thus, for each $U \in \U$ there exists a 
$y \in C$ such that $\Cy \subset U(y)$. The family $\C = \{\clCy\colon y \in C \}$ of closed 
sets has the finite intersection property and $\C$ contains small sets. Therefore 
\cite[Theorem 23, p. 193]{Ke}, there exists a unique $x \in \bigcap \C$. From the fact that 
$x \in \clCy$ and the lower semicontinuity of $q(\cdot,y)$ it follows that 
$q(x,y)=0$, $\yCx$. On the other hand, $x \in  \{z \in X\colon \psi(z) \leq \psi(y)\}$, for each 
$\yCx$, and consequently, $x \preq y$, $\yCx$, i.e. $x$ is the unique smallest element of $C$ (see 
(\ref{con2})). $\quad \square$

\par For a mapping $\psiX$ let us adopt
\begin{equation}
\label{con4}
 \beta=\inf\{\psi(z)\colon z \in X\} \mbox{, } B=\{z \in X\colon \psi(z) = \beta\}.
\end{equation}
Let us prove our basic result, which extends \cite[Theorem 21]{Pv}.
\begin{theorem}
\label{Th5}
Let $\XU$ be a uniform space, $q$ a $\U$-istance in $X$, while $q$ satisfies:
\begin{equation}
\label{con5}
 q(x,y) = 0 \mbox{ yields } x = y,  \quad \xyX
\end{equation}
or $\XU$ is Hausdorff. Assume that $\psiX$ is a lower semicontinuous mapping bounded below satisfying:
\begin{equation}
\label{con6}
 \psial=\{z \in X\colon \psi(z) \leq \alpha\} \mbox{ is complete for each } \alpha \in \mR,
\end{equation}
and that the following holds (see (\ref{con3}), (\ref{con4})):
\begin{eqnarray}
\label{con7}
 \mbox{for each } x \in X \setminus B  \mbox{ there exists a } \yXx \mbox{ such that } y \preq x .
\end{eqnarray}
Then for any $\xo \in X \setminus B$, each maximal chain $A \subset X$ containing $\xo$ 
(such a chain exists) has a unique smallest element $x$ which in addition satisfies:
\begin{enumerate}[(i)]
\item $\psi(x) = \inf \{\psi(z)\colon z \in X\}$ (i.e. $x \in B$),
\item $\psi(x) + q(x,\xo) - \psi(\xo) = \inf \{\psi(z) + q(z,\xo) - \psi(\xo)\colon z \in A \} \leq 0$,
\item $0 < \psi(y) + q(y,x) - \psi(x)$, $\quad \yXx$.
\end{enumerate}
\end{theorem}
\par Proof.
The relation $\preq$ is transitive and in view of Kuratowski's lemma \cite[p. 33]{Ke} for any 
$\xo \in X \setminus B$ (see (\ref{con7})) there exists a maximal chain $A$ containing $\xo$. Let us adopt 
$\alpha = \inf \{ \psi(z)\colon z \in A \}$ and suppose $\alpha < \psi(x)$, for each $x \in A$. 
Then there exists a sequence $\xnn$ in $A$ such that $\psixnn$ decreases to $\alpha$. Condition 
$\xm \preq \xn$ for $m < n$ would mean that

\[
 0 \leq q(\xm,\xn) \leq \psi(\xn) - \psi(\xm) < 0,
\]
which is impossible. Therefore, $\xn \preq \xm$ must be true ($\xn, \xm \in A$) and we obtain
\[
 q(\xn,\xm) \leq \psi(\xm) - \psi(\xn), \quad m,n \in \mN \mbox{, } m < n,
\]
i.e. $\xnn$ is a Cauchy sequence in $\Xq$, as $\psi$ is bounded below. Consequently, (see (\ref{con1c})) 
$\xnn$ is a Cauchy sequence in $\XU$. On the other hand, all $\xn$ for large 
$n$ are contained in $\{z \in X\colon \psi(z) \leq \psi(\xm) \}$ which is complete
(see (\ref{con6})). Therefore, $\xnn$ is convergent; let us say that $v \in \li \xn$. From the lower 
semicontinuity of $\psi$, $\qxm$ it follows that $\psi(v) \leq \alpha$ and 
\[
 \psi(v) + q(v,\xm) - \psi(\xm) \leq 0, \quad m \in \mN.
\]
On the other hand, for any $y \in A$ there exists an $\xn$ such that $\psi(\xn) < \psi(y)$ and
consequently, ($\xn, y \in A$)
\[
 \psi(\xn) + q(\xn,y) - \psi(y) \leq 0
\]
must hold, as
\[
 0 \leq q(y,\xn) \leq \psi(\xn) - \psi(y) < 0,
\]
is false. Now, we obtain (see (\ref{con1a}))
\[
 \psi(v) + q(v,y) - \psi(y) \leq \psi(v) + q(v,\xn) - \psi(\xn) + \psi(\xn) + q(\xn,y) - \psi(y) \leq 0. 
\]
This reasoning proves that there exists a $v \in A$ such that $\psi(v) = \alpha$. If $v$ is not a 
(unique) smallest element of $A$, then from $y \preq v$ for a $y \in A \setminus \{v\}$ it 
follows that
\[
 0 \leq q(y,v) \leq \psi(v) - \psi(y) \mbox{,}
\]
and we obtain $\psi(y) = \alpha $, $q(y,v) = 0$. Consequently, if  (\ref{con5}) is satisfied, then 
$x = v$ is the unique smallest element of $A$. Assume that $\XU$ is Hausdorff. Let us consider 
$C = \{z \in A\colon   q(z,v) = 0\}$. It is a maximal chain satisfying the assumptions of Lemma \ref{Le4} 
for $\psial$ in place of $X$ (see (\ref{con6})). The smallest element $x$ of $C$ is the unique smallest 
element of the maximal chain $A$. Conditions (ii), (iii) follow. Now, (\ref{con7}) yields $x \in B$ 
(otherwise $x$ would have another predecessor), i.e. (i) is satisfied. $\quad \square$
\par If $X$ is a complete uniform space, then condition (\ref{con6}) is satisfied and thus, 
\cite[Theorem 21]{Pv} is a consequence of Theorem \ref{Th5}.
\par Let us prove the following:
\begin{theorem}
\label{Th6}
Let $\XU$ be a uniform space, $q$ a $\U$-istance in $X$, while $q$ satisfies (\ref{con5}) or $\XU$ is 
Hausdorff. Assume that a lower semicontinuous mapping $\psiX$ is bounded 
below and  condition (\ref{con6}) holds. If a maximal chain $A \subset X$ (for $\preq$ as in 
condition (\ref{con3})) contains a point $\xo$, then $A$ has a unique smallest element $x$ which in addition 
satisfies: 
\begin{enumerate}[(i)]
\item $\psi(x) = \inf \{\psi(z)\colon z \in A\}$,
\item $\psi(x) + q(x,\xo) - \psi(\xo) = \inf \{\psi(z) + q(z,\xo) - \psi(\xo)\colon z \in A \} \leq 0$,
\item $0 < \psi(y) + q(y,x) - \psi(x)$, $\quad \yXx$,
\item $q(x,x)=0$, if $\ x \preq x$.
\end{enumerate}
\end{theorem}
\par Proof.
For any maximal chain $A$ containing $\xo$ we follow the proof of Theorem \ref{Th5}. If $x$ is the smallest 
element of $A$ and precedes itself, then the following holds
\[
 0 \leq q(x,x) = \psi(x) + q(x,x) - \psi(x) \leq 0.
\]
The last sentence of the proof of Theorem \ref{Th5} cannot be added, and we do not know if $x \in B$. $\quad \square$
\par The theorem of  Br\o ndsted \cite[Theorem 2]{Bro} was extended in the following 
way \cite[Theorem 23]{Pt} (we present its shorter formulation here, and adapted to the needs of the present paper):
\begin{theorem}
\label{Th7}
Let $\XU$ be a uniform space and $q\colon X \times X \to [0,\infty)$ a mapping satisfying:
\begin{enumerate}[(i)]
\item $q(z,x) \leq q(z,y) + q(y,x)$, $\quad \xyzX$,
\item $q(x,y) = 0$ yields $x = y$, $\quad \xyX$ (Br\o ndsted assumes  equivalence)
\item for each $U \in \mathcal{U}$ there exists a $\delta > 0$ such that 
 $q^{-1}([0,\delta)) \subset U$,
\item $q(\cdot,x)$ is lower semicontinuous, $\quad x \in X$.
\end{enumerate}
Assume that a lower semicontinuous mapping $\psiX$ is bounded below and  
$\{z \in X\colon \psi(z) \leq \alpha\}$ is complete for each $\alpha \in \mR$. Then for any 
$\xo \in X$ there exists an $x \in X$ such that 
$\psi(x) + q(x,\xo) - \psi(\xo) \leq q(\xo,\xo)$, and $0 < \psi(y) + q(y,x) - \psi(x)$, $\yXx$.
\end{theorem}
\par Let us compare Theorems \ref{Th7} and \ref{Th6}. If $\xo$ has no predecessor in $X$, then Theorem 
\ref{Th7} is trivial ($x=\xo$). Assume that $\xo$ has a predecessor. Conditions (i), (ii), (iv) are equivalent 
to (\ref{con1a}), (\ref{con5}), (\ref{con1b}), respectively. Condition (iii) is more restrictive than 
(\ref{con1c}) (e.g. for (iii) $\qxnxm < \epsilon$ is equivalent to $q(\xm,\xn) < \epsilon$). 
Consequently, Theorem \ref{Th6} is more general than \cite[Theorem 23]{Pt} and the Br\o ndsted theorem 
(the latter deals only with Hausdorff spaces).
\par Let us compare Theorem \ref{Th6} and \cite[Theorem 5]{Va}. V\'alyi considers mappings $d, \psi$ into 
a topological vector space ordered by the closed cone. In particular, for $d(x,y)=q(y,x)$, $\xyX$ and 
$\varphi = \psi$ our assumptions (\ref{con1a}), (\ref{con1b}) are the same as (ii), (iii) \cite[p. 30]{Va}; 
condition (\ref{con1c}) and V\'alyi's (5.3) are comparable (convergence in uniform space), and the requirements on 
$\varphi$ are the same as for our $\psi$. V\'alyi demands $\XU$ to be Hausdorff and $d(x,y) = 0$ iff $x=y$, 
while we assume that $q(x,y) = 0$ yields $x=y$ and only for the non-Hausdorff case. In place of the 
completeness of $\XU$ we use weaker condition (\ref{con6}). V\'alyi demands the transitive order relation to be 
also reflexive. Moreover, Theorem 5 in \cite{Va} states only our (iii) and a much simpler version of (ii). Thus 
V\'alyi's theorem in its ''real-valued'' version is weaker than Theorem \ref{Th6}.
\par Variational principles are related to fixed point theorems.
\par Let $2^Y$ be the family of all subsets of $Y$. We say that $\FXY$ is a (multivalued)
mapping if $F(x) \neq \emptyset$, for all $x \in X \neq \emptyset$.
\par The subsequent theorem extends the theorems of Caristi \cite[Theorem (2.1)']{Ca}, Takahashi 
\cite[Theorem 5]{Gr} and \cite[Theorem 24]{Pv}.
\begin{theorem}
\label{Th8}
Let $\XU$ be a uniform space, $q$ a $\U$-istance in $X$, with $q$ satisfying (\ref{con5}) or $\XU$ being 
Hausdorff. Assume that a lower semicontinuous mapping $\psiX$ is bounded 
below, condition (\ref{con6}) holds, $X \subset Y$ and $\FXY$ is a mapping satisfying 
(see (\ref{con3})):
\begin{eqnarray}
\label{con8}
 \mbox{for each } x \in X \mbox{ there exists a } y \in F(x) \mbox{ such that } y \preq x.
\end{eqnarray}
Then for any $\xo \in X$, each maximal chain $A \subset X$ containing $\xo$ (such a chain exists) 
has a unique smallest element $x$ which in addition satisfies conditions (i),(ii),(iii) of Theorem \ref{Th6}, 
$q(x,x) = 0$ and $x \in F(x)$. 
\end{theorem}
\par Proof.
In view of (\ref{con8}), Theorem \ref{Th6} here applies to nonempty chains. Now, condition (\ref{con8}) 
and (iii) mean $x \in F(x)$, and $q(x,x) = 0$. $\quad \square$
\par Under weaker assumptions we get the following extension of \cite[Theorem 23]{Pv}:
\begin{theorem}
\label{Th9}
Let $\XU$ be a uniform space, $q$ a $\U$-istance in $X$, with $q$ satisfying (\ref{con5}) or $\XU$ being 
Hausdorff. Assume that a lower semicontinuous mapping $\psiX$ is bounded 
below, condition (\ref{con6}) holds, $X \subset Y$, and $\FXY$ is a mapping satisfying 
(see (\ref{con3})):
\begin{eqnarray}
\label{con9}
 \mbox{for each } \xFx \mbox{ there exists a } \yXx  \mbox{ such that } y \preq x.
\end{eqnarray}
Then $F$ has a fixed point.
\end{theorem}
\par Proof.
Suppose $F$ has no fixed point. Then Theorem \ref{Th6} applies ($A$ is nonempty) and (iii) contradicts 
(\ref{con9}). $\quad \square$
\par Another consequence of Theorem \ref{Th6} is the following extension of Ekeland's variational principle 
\cite[Theorem 1]{Ek} and of more general results: \cite[Theorem 3]{Ka} and of \cite[Theorem 25]{Pv}.
\begin{theorem}
\label{Th10}
Let $\XU$ be a uniform space, $q$ a $\U$-istance in $X$, with $q$ satisfying (\ref{con5}) or $\XU$ being 
Hausdorff. Assume that a lower semicontinuous mapping $\psiX$ is bounded 
below and condition (\ref{con6}) holds. Then the following are satisfied:
\begin{enumerate}[(i)]
\item for each $\xo \in X$ there exists an $x \in X$ such that  
$\psi(x) \leq \psi(\xo)$ and $\psi(x) - q(y,x) < \psi(y)$, $\quad \yXx$, 
\item for any $\epsilon > 0$ and each $\xo \in X$ with $q(\xo,\xo) = 0$ and 
$\psi(\xo) \leq \epsilon + \inf\{\psi(z)\colon z \in X \}$ there exists an $x \in X$ such that 
$\psi(x) \leq \psi(\xo)$, $q(x,\xo) \leq 1$ and  
$\psi(x) - \epsilon q(y,x) < \psi(y)$, $\quad \yXx$.
\end{enumerate}
\end{theorem}
\par Proof.
The reasoning based on our Theorem \ref{Th6} is similar to the one presented in the proof of 
\cite[Theorem 3]{Ka}. The set $Y = \{z \in X\colon \psi(z) \leq \psi(\xo)\}$ is complete. Suppose that 
for each $x \in Y$ there exists a $y \in Y \setminus \{x\}$ such that 
$\psi(y) + q(y,x) - \psi(x) \leq 0$. Then by Theorem \ref{Th6} (iii) there exists an $x \in Y$ 
such that $0 < \psi(y) + q(y,x) - \psi(x)$ - a contradiction, i.e. (i) is proved for $Y$ in place of 
$X$. Condition (i) is trivial for $y \in X \setminus Y$.
\par Now, let us consider $\epsilon q$ in place of $q$. From $\xo \preq \xo$ and Theorem \ref{Th6} 
it follows that there exists an $x \in X$ such that 
$\psi(x) + \epsilon q(x,\xo) - \psi(\xo) \leq 0$ and 
$0 < \psi(y) + \epsilon q(y,x) - \psi(x)$, $\yXx$. From 
$\epsilon q(x,\xo) \leq \psi(\xo) - \psi(x)$ and the second assumption of (ii) we obtain
\[
 q(x,\xo) \leq [\psi(\xo) - \psi(x)]/\epsilon \leq [\psi(\xo) - \inf \{\psi(z)\colon z \in X\}]/\epsilon \leq 1.
 \quad \square
\]

\section*{Acknowledgements}
This work was partially supported by the Faculty of Applied Mathematics AGH UST statutory tasks within subsidy 
of the Polish Ministry of Science and Higher Education, grant no. 16.16.420.054.
\par 
\vspace{.1in}
\mbox{Faculty of Applied Mathematics} \linebreak
\mbox{AGH University of Science and Technology} \linebreak
\mbox{Al. Mickiewicza 30} \linebreak
\mbox{30-059 KRAK\'OW, POLAND} \linebreak
\mbox{E-mail: pasicki@agh.edu.pl}

\end{document}